\documentclass[11pt]{article}
\usepackage{tikz}
\usepackage[all]{xy}
\usepackage{amsfonts}
\usepackage{amsfonts}
\usepackage{amsfonts}
\usepackage{mathrsfs}
\usepackage{bm}
\usepackage{cite}
\usepackage[colorlinks,linkcolor=blue,citecolor=blue]{hyperref}
\usepackage{amssymb,amsmath} 
\usepackage{makecell}
\usepackage{tikz}
 \allowdisplaybreaks
\usetikzlibrary{arrows,shapes,chains}
 \allowdisplaybreaks

\catcode`!=11
\let\!int\int \def\int{\displaystyle\!int}
\let\!lim\lim \def\lim{\displaystyle\!lim}
\let\!sum\sum \def\sum{\displaystyle\!sum}
\let\!sup\sup \def\sup{\displaystyle\!sup}
\let\!inf\inf \def\inf{\displaystyle\!inf}
\let\!cap\cap \def\cap{\displaystyle\!cap}
\let\!max\max \def\max{\displaystyle\!max}
\let\!min\min \def\min{\displaystyle\!min}
\let\!frac\frac \def\frac{\displaystyle\!frac}
\catcode`!=12

\let\oldsection\section
\renewcommand\section{\setcounter{equation}{0}\oldsection}

\allowdisplaybreaks
\def\pf{\it{Proof.}\rm\quad}

\def\N{\mathbb{N}}

\def\su{\sum\limits_{n=1}^\infty}

\def\t{\widetilde{t}}
\def\S{\widetilde{S}}
\def\z{\zeta}
\def\T{{\bar T}}

\def\tt{\left(\frac{1-t}{1+t} \right)}

\newtheorem{defn}{Definition}[section]
\newtheorem{thm}{Theorem}[section]
\newtheorem{lem}[thm]{Lemma}
\newtheorem{cor}[thm]{Corollary}

\newtheorem{re}[thm]{Remark}

\begin{document}
\title {\bf Duality of Weighted Sum Formulas of Alternating Multiple $T$-Values}
\author{
{Ce Xu\thanks{Email: 19020170155420@stu.xmu.edu.cn}}\\[1mm]
\small $*$ School of Mathematics and Statistics, Anhui Normal University,\\ \small Wuhu 241000, P.R. China\\
[5mm]
Dedicated to professor Masanobu Kaneko on the occasion of his 60th birthday}

\date{}
\maketitle \noindent{\bf Abstract} Recently, a new kind of multiple zeta value level two $T({\bf k})$ (which is called multiple $T$-values) was introduced and studied by Kaneko and Tsumura. In this paper, we define a kind of alternating version of multiple $T$-values, and study several duality formulas of weighted sum formulas about alternating multiple $T$-values by using the methods of iterated integral representations and series representations. Some special values of alternating multiple $T$-values can also be obtained.
\\[2mm]


\section{Introduction}

We begin with some basic notations. A finite sequence ${\bf k} := (k_1,\ldots, k_r)$ of positive integers is called an index. We
put
\[|{\bf k}|:=k_1+\cdots+k_r,\quad d({\bf k}):=r,\]
and call them the weight and the depth of ${\bf k}$, respectively. If $k_r>1$, ${\bf k}$ is called admissible. Let $I(k,r)$ be the set of all indices of weight $k$ and depth $r$.

For an admissible index ${\bf k}=(k_1,\ldots,k_r)$, the multiple zeta values (abbr. MZVs) are defined by
\begin{align}
\zeta(k_1,k_2,\ldots,k_r):=\sum\limits_{0<m_1<\cdots<m_r } \frac{1}{m_1^{k_1}m_2^{k_2}\cdots m_r^{k_r}}.
\end{align}
When $r=1$ and $k_1=k$, $\z(k)$ stands for Riemann zeta value, which is defined by
\begin{align*}
\z(k):=\su \frac{1}{n^k} \quad (k\in \N\setminus\{1\}).
\end{align*}
The systematic study of multiple zeta values began in the early 1990s with the works of Hoffman \cite{H1992} and Zagier \cite{DZ1994}. After that it has been attracted a lot of research on them in the last three decades (see, for example, the book of Zhao \cite{Z2016}).

So far, surprisingly little work has been done on several variants of multiple zeta values. In recent papers \cite{KTA2018,KTA2019}, Kaneko and Tsumura introduced and studied a new kind of multiple zeta values of level two
\begin{align}
T(k_1,k_2,\ldots,k_r):&=2^r \sum_{0<m_1<\cdots<m_r\atop m_j\equiv j\ {\rm mod}\ 2} \frac{1}{m_1^{k_1}m_2^{k_2}\cdots m_r^{k_r}}\nonumber\\
&=2^r\sum\limits_{0<n_1<\cdots<n_r} \frac{1}{(2n_1-1)^{k_1}(2n_2-2)^{k_2}\cdots (2n_r-r)^{k_r}}\label{a2}
\end{align}
for an admissible index ${\bf k}=(k_1,\ldots,k_r)$, which is called the `multiple $T$-value' (MTV for short). This is in contrast to Hoffman's multiple $t$-value (MtV for short) defined by
\begin{align}
t(k_1,k_2,\ldots,k_r):&=\sum_{0<m_1<\cdots<m_r\atop \forall m_j:\ {\rm odd}} \frac{1}{m_1^{k_1}m_2^{k_2}\cdots m_r^{k_r}},
\end{align}
which was introduced and studied in his recent paper \cite{H2019} as another variant of multiple zeta values of level two.

In (\ref{a2}), we put a bar
on top of $k_j$ if there is a sign $(-1)^{n_j}$ appearing in the denominator on the right. For example,
\begin{align*}
T(k_1,{\overline {k_2}},k_3,{\overline {k_4}}):=2^4\sum\limits_{0<n_1<n_2<n_3<n_4} \frac{(-1)^{n_2+n_4}}{(2n_1-1)^{k_1}(2n_2-2)^{k_2}(2n_3-3)^{k_4}(2n_4-4)^{k_4}}.
\end{align*}
The sums of types above (one of more the $k_j$ barred) are called the alternating multiple $T$-values.
In particular, we let
\begin{align}
\T(k_1,k_2,\ldots,k_r):=2^r\sum\limits_{0<n_1<\cdots<n_r} \frac{(-1)^{n_r}}{(2n_1-1)^{k_1}(2n_2-2)^{k_2}\cdots (2n_r-r)^{k_r}},
\end{align}
for an index ${\bf k}=(k_1,\ldots,k_r)$. In \cite[Corollary 3.4]{WX2020}, the author and Wang proved that the four alternating linear sums
\begin{align*}
&(1-(-1)^{p+q})\su \frac{h^{(p)}_n}{n^q},\quad (1+(-1)^{p+q})\su \frac{{\bar h}^{(p)}_n}{n^q},\\
&(1+(-1)^{p+q})\su \frac{{\bar h}^{(p)}_n}{n^q}(-1)^{n-1},\quad (1-(-1)^{p+q})\su \frac{h^{(p)}_n}{n^q}(-1)^{n-1}
\end{align*}
can be expressed in terms of (alternating) Riemann zeta values and (alternating) $t$-values.
Here $h_n^{(p)}$ stands for odd harmonic number of order $p$ defined by
\begin{align*}
h_n^{(p)}:=\sum\limits_{k=1}^n \frac{1}{(k-1/2)^p},\quad h_n\equiv h_n^{(1)} \quad{\rm and}\quad h_0^{(p)}:=0
\end{align*}
and ${\bar h}_n^{(p)}$ stands for alternating odd harmonic number of order $p$ defined by
\begin{align*}
{ \bar h}_n^{(p)}: = \sum\limits_{k=1}^n {\frac{(-1)^{k-1}}{{{(k-1/2)^p}}}},\quad {\bar h}_n\equiv{\bar h}^{(1)}_n,\quad {\bar h}_0^{(p)}:=0.
\end{align*}
Hence, according to definitions, we know that the two (alternating) double $T$-values
\[T(k_1,k_2)=\frac{1}{2^{k_1+k_2-2}}\su \frac{h_n^{(k_1)}}{n^{k_2}},\quad T(k_1,\overline{k_2})=\frac{1}{2^{k_1+k_2-2}}\su \frac{h_n^{(k_1)}}{n^{k_2}}(-1)^{n-1}\]
can be evaluated by (alternating) Riemann zeta values and (alternating) $t$-values if weight $k_1+k_2$ is odd (Note that the explicit evaluation of $T(k_1,k_2)$ with weight $k_1+k_2$ odd was first found by Kaneko-Tsumura in \cite{KTA2018}). The two (alternating) double $T$-values
\[T(\overline{k_1},k_2)=-\frac{1}{2^{k_1+k_2-2}}\su \frac{{\bar h}_n^{(k_1)}}{n^{k_2}},\quad T(\overline{k_1},\overline{k_2})=-\frac{1}{2^{k_1+k_2-2}}\su \frac{{\bar h}_n^{(k_1)}}{n^{k_2}}(-1)^{n-1}\]
can be evaluated by (alternating) Riemann zeta values and (alternating) $t$-values if weight $k_1+k_2$ is even.
Here the alternating Riemann zeta value ${\bar \z}(k)$ is denoted by
\begin{align}
{\bar \z}(k):=\su \frac{(-1)^{n-1}}{n^k} \quad \text{and} \quad {\bar \z}(0):=\frac{1}{2}\quad (k\in \N),
\end{align}
and the alternating $t$-value ${\bar t}(k)$ is denoted by
\begin{align}
{\bar t}(k):=\su \frac{(-1)^{n-1}}{(n-1/2)^k}\quad (k\in\N).
\end{align}
Clearly, ${\bar \z}(k)=(1-2^{1-k})\z(k)$ if $k>1$ and $\T(k)=T({\bar k})=-{\bar t}(k)/{2^{k-1}}$. Note that from \cite{G1989}, for nonnegative integer $k$, we have the generating function of ${\bar t}(2k+1)$,
\begin{align*}
\frac{\pi}{\cos(\pi s)}=2\sum\limits_{k=0}^{\infty} {\bar t}(2k+1)s^{2k}=\sum\limits_{k=0}^{\infty} \frac{(-1)^kE_{2k}\pi^{2k+1}}{(2k)!}s^{2k},
\end{align*}
where $E_{2k}$ is Euler number. Thus, we compute
\begin{align}
{\bar t}(2k+1)=\frac{(-1)^kE_{2k}\pi^{2k+1}}{2(2k)!}\quad (k\geq 0).
\end{align}

For multiple zeta values, the classical sum formula is widely known and its variants are enormous (see \cite[Chapter 5]{Z2016}), for example, Zagier \cite{DZ1994} conjectured the following weighted sum formula for MZVs that for an admissible index ${\bf k}$,
\begin{align*}
\sum_{{\bf k}\in I(k,r)}\z(k_1,k_2,\ldots,k_r)=\z(k)\quad (k\geq 2)
\end{align*}
which was first proved by Granville in \cite{G1997}. Some other related results for weighted sums of MZVs may be seen in the works of \cite{GKZ2006,G2016,GL2015,H2017,LQ2019,M2013,N2009} and references therein. For Hoffman's multiple $t$-values, many mathematicians studied the weighed sums formulas problem of MtVs, for example see
\cite{LX2020,SC2012,SJ2017,Z2015}.
In contrast, we know very little about the weighed sum formula for MTVs. For MTVs, only certain formulas are found in depths 2 and 3, for detail see Kaneko and Tsumura \cite[Theorems 3.2 and 3.3]{KTA2019}. Moreover, in \cite{KTA2019}, Kaneko and Tsumura also conjectured an analogue of Machide's formula for MTVs, see \cite[Conjecture 4.6]{KTA2019}. In this paper, we will discuss a kind of weighted sums of alternating MTVs.
For an index ${\bf k}=(k_1,\ldots,k_r)$ and positive integer $m$, we define a weighted sum of alternating MTVs with arguments of weight $k$ and depth $r$ by
\begin{align}
W_m(k,r):=\sum_{{\bf k}\in I(k,r)}\binom{k_r+m-2}{m-1} \T(k_1,\ldots,k_{r-1},k_r+m-1),
\end{align}
if $m=1$, we let $W(k,r):=W_1(k,r)$. Obviously, with the help of \cite[Eq.(3.10)]{WX2020}, by a direct calculation, we arrive at the conclusion that if $k>1$ is an odd, then the weighted sum
$W(k,2)$
can be expressed in terms of (alternating) Riemann zeta values (for detailed, see Theorem \ref{thmc4} in the present paper).

For positive integers $j$ and $p$ with $1\leq j\leq p$, we let
\begin{align}
Z(j,p):=\sum_{k=1}^{p-j}(-2)^k\sum\limits_{i_0<i_1<\cdots<i_{k-1}<i_k,\atop i_0=j,i_k=p} \prod\limits_{l=1}^k {\bar \z}(2i_l-2i_{l-1}),\quad Z(p,p):=1.
\end{align}
It is clear that
$Z(j,p)=a_{j,p}\zeta(2p-2j)\quad (a_{j,p}\in \mathbb{Q})$.
For instance, we have $Z(j,j+1)=-\z(2),\ Z(j,j+2)=3\z(4)/4$ and $Z(j,j+3)=-3\z(6)/16$. In particular, the author and Zhao \cite{XZ2020} shown that
\begin{align}
Z(j,p)=\frac{(-1)^{p-j}\pi^{2p-2j}}{(2p-2j+1)!}.
\end{align}

The purpose of present paper is to evaluate the duality of products of weighted sum $W(k,r)$ and $Z(j,p)$ by using the method of iterated integrals. Our main results are following three theorems.

\begin{thm}\label{thm1} For positive integers $m$ and $p$, we have
\begin{align}\label{a9}
(-1)^m \sum_{j=1}^p W(2j+2m-1,2m)Z(j,p)=(-1)^p \sum_{j=1}^m W(2j+2p-1,2p)Z(j,m).
\end{align}
\end{thm}

\begin{thm}\label{thm2} For positive integers $m$ and $p$, we have
\begin{align}\label{a10}
(-1)^m \sum_{j=1}^p W(2j+2m-3,2m-1)Z(j,p)=(-1)^p \sum_{j=1}^m W(2j+2p-3,2p-1)Z(j,m).
\end{align}
\end{thm}

\begin{thm}\label{thm3} For positive integers $k,m$ and $r$, we have
\begin{align}\label{a11}
&W_m(k+r-1,r)+(-1)^r W_k(m+r-1,r)\nonumber\\
&=\sum_{j=1}^{r-1} (-1)^{j-1} T(\{1\}_{k-2},{\bar 1},\{1\}_{j-1},{\bar 1})T(\{1\}_{m-2},{\bar 1},\{1\}_{r-j-1},{\bar 1}),
\end{align}
where $(\{1\}_{-1},{\bar 1}):=\emptyset$, and $\{l\}_m$ denotes the sequence $\underbrace{l,\ldots,l}_{m \text{\;times}}$.
\end{thm}

Clearly, from Theorems \ref{thm1}, \ref{thm2} and \ref{thm3}, many explicit relations about weighted sums of alternating MTVs can be established. For example, by straightforward calculations, we obtain
\begin{align*}
W(5,4)+W(5,2)=\frac{7}{4}\z(2)\z(3)\quad\text{and}\quad
W(7,5)=\frac{\pi^7}{7680}+\z(2)W(5,3).
\end{align*}

\section{Preliminaries}

In this section we will establish several fundamental formulas for weighed sums $W(k,r)$.

Kaneko and Tsumura \cite{KTA2018,KTA2019} introduced the following a new kind of multiple polylogarithm function
\begin{align}\label{b1}
&{\rm A}(k_1,k_2,\ldots,k_r;z): = 2^r\sum\limits_{1 \le {n_1} <  \cdots  < {n_r}\atop n_j\equiv j\ {\rm mod}\ 2} {\frac{{{z^{{n_r}}}}}{{n_1^{{k_1}}n_2^{{k_2}} \cdots n_r^{{k_r}}}}}\quad (|z|<1),
\end{align}
where $k_1,\ldots,k_r$ are positive integers. If $k_r>1$, we allow $z=1$\ (In \cite{KTA2018}, $2^{-r}{\rm A}(k_1,\ldots,k_r;z)$ is denoted by ${\rm Ath}(k_1,\ldots,k_r;z)$.). When $r=1$, the polylogarithm of level 2 ${\rm A}(k;z)$ was considered by Sasaki in \cite{S2012}.

From definition, it is easy to find that
\begin{align}\label{b2}
\frac{d}{dz}{\mathrm{A}}({{k_1}, \cdots ,k_{r-1},{k_r}}; z)= \left\{ {\begin{array}{*{20}{c}} \frac{1}{z} {\mathrm{A}}({{k_1}, \cdots ,{k_{r-1}},{k_r-1}};z)
   {\ \ (k_r\geq 2),}  \\
   {\frac{2}{1-z^2}{\mathrm{A}}({{k_1}, \cdots ,{k_{r-1}}};z)\;\;\;\ \ \ (k_r = 1).}  \\
\end{array} } \right.
\end{align}
Hence, by (\ref{b2}), we deduce the iterated integral expression
\begin{align}\label{b3}
&{\mathrm{A}}({{k_1}, \cdots,k_{r-1} ,{k_r}};z)=\int\limits_{0}^z \underbrace{\frac{dt}{t}\cdots\frac{dt}{t}}_{k_r-1}\frac{2dt}{1-t^2}\underbrace{\frac{dt}{t}\cdots\frac{dt}{t}}_{k_{r-1}-1}\frac{2dt}{1-t^2}\cdots
\underbrace{\frac{dt}{t}\cdots\frac{dt}{t}}_{k_1-1}\frac{2dt}{1-t^2}\nonumber\\
&=2^r\left\{\prod\limits_{j=1}^r\frac{(-1)^{k_j-1}}{\Gamma(k_j)}\right\}\int\nolimits_{D_r(z)} \frac{\log^{k_1-1}\left(\frac{t_1}{t_2}\right)\cdots \log^{k_{r-1}-1}\left(\frac{t_{r-1}}{t_r}\right)\log^{k_r-1}\left(\frac{t_r}{z}\right)}{(1-t_1^2)\cdots (1-t_{r-1}^2)(1-t_r^2)}dt_1\cdots dt_r,
\end{align}
where $D_{r}(z):=\{(t_1,\ldots,t_r)\mid 0<t_1<\cdots<t_r<z\}\quad (r\in\N).$ Here,
\[\int\limits_{0}^z f_k(t)dtf_{k-1}(t)dt \cdots f_1(t_1)dt:=\int\limits_{D_k(z)} f_k(t_k)f_{k-1}(t_{k-1}) \cdots f_1(t_1)dt_1\cdots dt_{k-1}dt_k.\]
Further, setting $t_j\rightarrow zt_j$ in (\ref{b3}) gives
\begin{align}\label{b4}
&{\rm A}(k_1,\ldots,k_{r-1},k_r;z)\nonumber\\
&=2^rz^r\left\{\prod\limits_{j=1}^r\frac{(-1)^{k_j-1}}{\Gamma(k_j)}\right\}\int\nolimits_{D_r(1)} \frac{\log^{k_1-1}\left(\frac{t_1}{t_2}\right)\cdots \log^{k_{r-1}-1}\left(\frac{t_{r-1}}{t_r}\right)\log^{k_r-1}\left(t_r\right)}{(1-z^2t_1^2)\cdots (1-z^2t_{r-1}^2)(1-z^2t_r^2)}dt_1\cdots dt_r.
\end{align}
Applying the changes of variables $t_j\rightarrow \frac{1-t_{r+1-j}}{1+t_{r+1-j}}$ to (\ref{b4}) yields
\begin{align}\label{b5}
&{\rm A}(k_1,\ldots,k_{r-1},k_r;z)\nonumber\\
&=4^rz^r\left\{\prod\limits_{j=1}^r\frac{(-1)^{k_j-1}}{\Gamma(k_j)}\right\}\int\nolimits_{D_r(1)} \frac{\log^{k_1-1}\left(\frac{(1-t_r)(1+t_{r-1})}{(1+t_r)(1-t_{r-1})} \right)}{(1+t_r)^2-z^2(1-t_r)^2}\cdots  \frac{\log^{k_{r-1}-1}\left(\frac{(1-t_2)(1+t_{1})}{(1+t_2)(1-t_{1})} \right)}{(1+t_2)^2-z^2(1-t_2)^2}\nonumber\\
&\quad\quad\quad\quad\quad\quad\quad\quad\quad\quad\quad\quad \times \frac{\log^{k_{r}-1}\left(\frac{1-t_1}{1+t_1} \right)}{(1+t_1)^2-z^2(1-t_1)^2}dt_1dt_2\cdots dt_r.
\end{align}
As is well-known, we can regard the function ${\rm A}(k_1,k_2,\ldots,k_r;z)$ as a single-valued holomorphic
function in the simply connected domain $\mathbb{C}\setminus \{(-\infty,-1]\cup [1,+\infty)\}$, via the process of iterated integration
starting with ${\rm A}(1;z)=\int_0^z 2dt/{(1-t^2)}$.

Next, we use (\ref{b5}) to establish four identities involving weighted sum $W(k,r)$, alternating Riemann zeta values and infinite series whose general $n$-terms is a product of $(-1)^n/n$, multiple $T$-harmonic sum and multiple $S$-harmonic sum. First, we need to give the definitions of multiple $T$-harmonic sum and multiple $S$-harmonic sum. For an index ${{\bf k}}:= (k_1,\ldots, k_r)$, let
\begin{align*}
{\bf k}_{2m-1}:=(k_1,k_2,\ldots,k_{2m-1})\quad \text{and}\quad {\bf k}_{2m}:=(k_1,k_2,\ldots,k_{2m}).
\end{align*}
For positive integers $n_1,n_2,\ldots,n_{r}$ and $n$, if $r=2m-1$ is an odd, we define
\begin{align*}
&D_n({{\bf n}_{2m-1}}):=\left\{(n_1,n_2,\ldots,n_{2m-1},n)\mid 0<n_1\leq n_2 <\cdots\leq n_{2m-2}<n_{2m-1}\leq n \right\},\ (n\geq m)\\
&E_n({{\bf n}_{2m-1}}):=\left\{(n_1,n_2,\ldots,n_{2m-1},n)\mid 1\leq n_1<n_2\leq \cdots< n_{2m-2}\leq n_{2m-1}< n \right\},\ (n>m)
\end{align*}
and if $r=2m$ is an even, we define
\begin{align*}
&D_n({{\bf n}_{2m}}):=\left\{(n_1,n_2,\ldots,n_{2m},n)\mid 0<n_1\leq n_2 <\cdots\leq n_{2m-2}<n_{2m-1}\leq n_{2m}<n \right\},\ (n> m)\\
&E_n({{\bf n}_{2m}}):=\left\{(n_1,n_2,\ldots,n_{2m},n)\mid 1\leq n_1<n_2\leq \cdots< n_{2m-2}\leq n_{2m-1}< n_{2m}\leq n \right\},\ (n>m).
\end{align*}

\begin{defn}(\cite{XZ2020})\label{lem1} For positive integer $m$, the multiple $T$-harmonic sums ({\rm MTHSs} for short) and multiple $S$-harmonic sums ({\rm MSHSs} for short) are defined by
\begin{align}
&T_n({{\bf k}_{2m-1}}):= \sum_{D_n({{\bf n}_{2m-1}})} \frac{2^{2m-1}}{(\prod_{j=1}^{m-1} (2n_{2j-1}-1)^{k_{2j-1}}(2n_{2j})^{k_{2j}})(2n_{2m-1}-1)^{k_{2m-1}}},\label{MOT}\\
&T_n({{\bf k}_{2m}}):= \sum_{D_n({{\bf n}_{2m}})} \frac{2^{2m}}{\prod_{j=1}^{m} (2n_{2j-1}-1)^{k_{2j-1}}(2n_{2j})^{k_{2j}}},\label{MET}\\
&S_n({{\bf k}_{2m-1}}):= \sum_{E_n({{\bf n}_{2m-1}})} \frac{2^{2m-1}}{(\prod_{j=1}^{m-1} (2n_{2j-1})^{k_{2j-1}}(2n_{2j}-1)^{k_{2j}})(2n_{2m-1})^{k_{2m-1}}},\label{MOS}\\
&S_n({{\bf k}_{2m}}):= \sum_{E_n({{\bf n}_{2m}})} \frac{2^{2m}}{\prod_{j=1}^{m} (2n_{2j-1})^{k_{2j-1}}(2n_{2j}-1)^{k_{2j}}},\label{MES}
\end{align}
where $T_n({{\bf k}_{2m-1}}):=0$ if $n<m$, and $T_n({{\bf k}_{2m}})=S_n({{\bf k}_{2m-1}})=S_n({{\bf k}_{2m}}):=0$ if $n\leq m$. Moreover, for convenience we let $T_n(\emptyset)=S_n(\emptyset):=1$. We call (\ref{MOT}) and (\ref{MET}) are multiple $T$-harmonic sums, and call (\ref{MOS}) and (\ref{MES}) are multiple $S$-harmonic sums.
\end{defn}

Clearly, according to the definitions of ${\rm A}(k_1,k_2,\ldots,k_r;z)$ and MTHSs, by an elementary calculation, we can find that
\begin{align}
&{\rm A}({\bf k}_{2m-1};z)=2\su \frac{T_n({\bf k}_{2m-2})}{(2n-1)^{k_{2m-1}}}z^{2n-1}
\end{align}
and
\begin{align}
&{\rm A}({\bf k}_{2m};z)=2\su \frac{T_n({\bf k}_{2m-1})}{(2n)^{k_{2m}}}z^{2n}.
\end{align}
From the definitions of alternating MTVs and ${\rm A}(k_1,k_2,\ldots,k_r;z)$, we have
\begin{align}\label{b12}
&{\rm A}({\bf k}_{2m-1};i)=i(-1)^m {\bar T}({\bf k}_{2m-1})
\end{align}
and
\begin{align}\label{b13}
&{\rm A}({\bf k}_{2m};i)=(-1)^m {\bar T}({\bf k}_{2m}),
\end{align}
where $i$ is imaginary unit. Kaneko and Tsumura \cite{KTA2018} gave the following relation
\begin{align}\label{b14}
{\rm A}(\{1\}_r;z)=\frac{1}{r!} \left({\rm A}(1;z)\right)^r=\frac{(-1)^r}{r!}\log^r\left(\frac{1-z}{1+z} \right).
\end{align}
It is clear that $A(1;i)=i\frac{\pi}{2}$. Hence, using (\ref{b12})-(\ref{b14}), a simple calculation yields
\begin{align}\label{b15}
{\bar T}(\{1\}_r)=\frac{(-1)^r}{r!}\left(\frac{\pi}{2}\right)^r.
\end{align}

Next, in order to state our main results, we need to give a lemma.
\begin{lem}(\cite[Theorem 2.1]{XZ2020})\label{lem1} For positive integers $m$ and $n$, the following identities hold
\begin{align}
&\begin{aligned}
\int_{0}^1 t^{2n-2} \log^{2m}\tt dt&= \frac{2(2m)!}{2n-1} \sum_{j=0}^m {\bar \z}(2m-2j)T_n(\{1\}_{2j}),\label{ee}
\end{aligned}\\
&\begin{aligned}
\int_{0}^1 t^{2n-2} \log^{2m-1}\tt dt&= -\frac{2(2m-1)!}{2n-1} \sum_{j=0}^{m-1} {\bar \z}(2m-1-2j)T_n(\{1\}_{2j})\\&\quad-\frac{(2m-1)!}{2n-1} S_n(\{1\}_{2m-1}),\label{eo}
\end{aligned}\\
&\begin{aligned}
\int_{0}^1 t^{2n-1} \log^{2m}\tt dt&=\frac{(2m)!}{n} \sum_{j=0}^{m-1} {\bar \z}(2m-1-2j)T_n(\{1\}_{2j+1})\\&\quad+\frac{(2m)!}{2n} S_n(\{1\}_{2m}),\label{oe}
\end{aligned}\\
&\begin{aligned}
\int_{0}^1 t^{2n-1} \log^{2m-1}\tt dt&= -\frac{(2m-1)!}{n} \sum_{j=0}^{m-1} {\bar \z}(2m-2-2j)T_n(\{1\}_{2j+1}),\label{oo}
\end{aligned}
\end{align}
where ${\bar \zeta}(0)$ should be interpreted as $1/2$ wherever it occurs.
\end{lem}

Using the iterated integral identity (\ref{b5}) with the help of Lemma \ref{lem1}, we can get the following theorem.
\begin{thm} For positive integers $m$ and $p$,
\begin{align}
&\begin{aligned}\label{b20}
W(2p+2m-1,2m)=2(-1)^m \sum_{j=0}^{p-1} {\bar \z}(2p-2-2j)\su \frac{T_n(\{1\}_{2m-1})T_n(\{1\}_{2j+1})}{n}(-1)^n,
\end{aligned}\\
&\begin{aligned}\label{b21}
W(2p+2m-2,2m)&=2(-1)^m \sum_{j=0}^{p-2} {\bar \z}(2p-3-2j)\su \frac{T_n(\{1\}_{2m-1})T_n(\{1\}_{2j+1})}{n}(-1)^n\\
&\quad+(-1)^m \su \frac{T_n(\{1\}_{2m-1})S_n(\{1\}_{2p-2})}{n}(-1)^n,
\end{aligned}\\
&\begin{aligned}\label{b22}
W(2p+2m-2,2m-1)&=4(-1)^{m-1} \sum_{j=0}^{p-1} {\bar \z}(2p-1-2j)\su \frac{T_n(\{1\}_{2m-2})T_n(\{1\}_{2j})}{2n-1}(-1)^n\\
&\quad+2(-1)^{m-1} \su \frac{T_n(\{1\}_{2m-2})S_n(\{1\}_{2p-1})}{2n-1}(-1)^n,
\end{aligned}\\
&\begin{aligned}\label{b23}
W(2p+2m-3,2m-1)=4(-1)^{m-1} \sum_{j=0}^{p-1} {\bar \z}(2p-2-2j)\su \frac{T_n(\{1\}_{2m-2})T_n(\{1\}_{2j})}{2n-1}(-1)^n.
\end{aligned}
\end{align}
\end{thm}
\pf By using (\ref{b12}) and (\ref{b13}) we know that
\begin{align*}
{\bar T}(k_1,k_2,\ldots,k_r)=i^r {\rm A}(k_1,k_2,\ldots,k_r;i).
\end{align*}
Then summing both sides of above identity and using (\ref{b5}), we have
\begin{align}\label{bb1}
W(k+r-1,r)=\frac{(-1)^{k+r-1}}{(k-1)!}2^r \int_{D_r(1)} \frac{\log^{k-1}\left(\frac{1-t_r}{1+t_r}\right)}{(1+t_1^2)\cdots (1+t_r^2)}dt_1\cdots dt_r.
\end{align}
Expand the $(1+t_j)^{-1}$ into geometric series to deduce that
\begin{align*}
W(k+r-1,r)=\frac{(-1)^{k-1}}{(k-1)!}2^r \sum_{0<m_1<\cdots<m_r} \frac{(-1)^{m_r} \int_0^1 t^{2m_r-r-1}\log^{k-1}\tt dt}{(2m_1-1)(2m_2-2)\cdots (2m_{r-1}-r+1)}.
\end{align*}
By straightforward calculations, if $r=2m$ and $k=2p$, then
\begin{align*}
W(2p+2m-1,2m)=\frac{2(-1)^{m-1}}{(2p-1)!} \su (-1)^n T_n(\{1\}_{2m-1}) \int_0^1 t^{2n-1}\log^{2p-1}\tt dt,
\end{align*}
if $r=2m$ and $k=2p-1$, then
\begin{align*}
W(2p+2m-2,2m)=\frac{2(-1)^{m}}{(2p-2)!} \su (-1)^n T_n(\{1\}_{2m-1}) \int_0^1 t^{2n-1}\log^{2p-2}\tt dt,
\end{align*}
if $r=2m-1$ and $k=2p$, then
\begin{align*}
W(2p+2m-2,2m-1)=\frac{2(-1)^{m}}{(2p-1)!} \su (-1)^n T_n(\{1\}_{2m-2}) \int_0^1 t^{2n-2}\log^{2p-1}\tt dt,
\end{align*}
if $r=2m-1$ and $k=2p-1$, then
\begin{align*}
W(2p+2m-3,2m-1)=\frac{2(-1)^{m-1}}{(2p-2)!} \su (-1)^n T_n(\{1\}_{2m-2}) \int_0^1 t^{2n-2}\log^{2p-2}\tt dt.
\end{align*}
Thus, with the help of formulas (\ref{ee})-(\ref{oo}), we may easily deduce these desired evaluations.\hfill$\square$

Next, we give several special cases. Setting $m=p=1$ in (\ref{b20}) yields
\begin{align}\label{b24}
\su \frac{T_n^2(1)}{n}(-1)^{n-1}={\T}(1,2)+\T(2,1)=\frac{7}{4}\z(3),
\end{align}
where we used these two well-known results
\begin{align*}
\T(1,2)=-\frac{7}{4}\z(3)+\frac{\pi}{4}{\bar t}(2)\quad \text{and}\quad \T(2,1)=\frac{7}{2}\z(3)-\frac{\pi}{4}{\bar t}(2).
\end{align*}
Letting $m=1$ in (\ref{b23}) gives
\begin{align}
\T(2p-1)=2\sum_{j=1}^{p-1} (-1)^j{\bar \z}(2p-2-2j) \T(\{1\}_{2j+1}).
\end{align}

Now, we end this section by a theorem.
\begin{thm}\label{thmb3} For positive integers $r$ and $k$,
\begin{align}
\T(\{1\}_{r-1},k)=\sum_{j=1}^r (-1)^{j-1} \T(\{1\}_{r-j})W(k+j-1,j).
\end{align}
\end{thm}
\pf According to the definitions of alternating MTVs and using (\ref{b5}) with $z=i,\ k_1=\cdots=k_{r-1}=1$ and $k_r=1$, we have
\begin{align}\label{bb2}
\T(\{1\}_{r-1},k)=\frac{(-1)^{k+r-1}}{(k-1)!}2^r \int_{D_r(1)} \frac{\log^{k-1}\left(\frac{1-t_1}{1+t_1}\right)}{(1+t_1^2)\cdots (1+t_r^2)}dt_1\cdots dt_r.
\end{align}
Note that the two integrals
\begin{align*}
\int_{0<t_1<\cdots<t_r<x} \frac{dt_1\cdots dt_r}{(1+t_1^2)\cdots(1+t_r^2)}=\frac{1}{r!} \arctan(x)^r
\end{align*}
and
\begin{align*}
\int_{x<t_1<\cdots<t_r<1} \frac{dt_1\cdots dt_r}{(1+t_1^2)\cdots(1+t_r^2)}=\frac{1}{r!} \left(\frac{\pi}{4}-\arctan(x)\right)^r,
\end{align*}
then the identities (\ref{bb1}) and (\ref{bb2}) can be rewritten in the forms
\begin{align*}
W(k+r-1,r)=\frac{(-1)^{k+r-1}2^r}{(k-1)!(r-1)!} \int_{0}^1 \frac{\log^{k-1}\tt \arctan(t)^{r-1}}{1+t^2}dt
\end{align*}
and
\begin{align*}
\T(\{1\}_{r-1},k)=\frac{(-1)^{k+r-1}2^r}{(k-1)!(r-1)!} \int_{0}^1 \frac{\log^{k-1}\tt \left(\frac{\pi}{4}- \arctan(t)\right)^{r-1}}{1+t^2}dt.
\end{align*}
Thus, combining the two equations above, by an elementary calculation, we prove the desired result.\hfill$\square$

\section{Proofs of Theorems \ref{thm1} and \ref{thm2}}

In this section, we give the proofs of Theorems \ref{thm1} and \ref{thm2}. We need the following a lemma.
\begin{lem}(\cite[Lemma 5.1]{XZ2020})\label{lemc1}  Let $A_{p,q}, B_p, C_p\ (p,q\in \N)$ be any complex sequences. If
\begin{align}\label{e1}
\sum\limits_{j=1}^p A_{j,p}B_j=C_p\quad\text{and}\quad A_{p,p}:=1,
\end{align}
hold, then
\begin{align}\label{e2}
B_p=\sum\limits_{j=1}^p C_j \sum\limits_{k=1}^{p-j} (-1)^k \left\{\sum\limits_{i_0<i_1<\cdots<i_{k-1}<i_k,\atop i_0=j,i_k=p} \prod\limits_{l=1}^k A_{i_{l-1},i_l}\right\},
\end{align}
holds, where $\sum\limits_{k=1}^0 (\cdot):=1$.
\end{lem}

Hence, from (\ref{b20}), (\ref{b23}) and Lemma \ref{lemc1}, we can get the following two theorems.
\begin{thm} For positive integers $m$ and $p$,
\begin{align}\label{c3}
\su \frac{T_n(\{1\}_{2m-1})T_n(\{1\}_{2p-1})}{n}(-1)^n=(-1)^m \sum_{j=1}^p W(2j+2m-1,2m)Z(j,p).
\end{align}
\end{thm}
\pf  Setting
\begin{align*}
&A_{j,p}=2{\bar \z}(2p-2j),\quad B_j=\su \frac{T_n(\{1\}_{2m-1})T_n(\{1\}_{2j-1})}{n}(-1)^n
\end{align*}
and
\begin{align*}
&C_p=(-1)^mW(2p+2m-1,2m)
\end{align*}
in Lemma \ref{lemc1} and using (\ref{b20}) gives the desired formula.\hfill$\square$

Putting $m=1$ and $p=2$ in (\ref{c3}), we get
\begin{align*}
\su \frac{T_n(1)T_n(1,1,1)}{n}(-1)^n&=\frac{7}{4}\z(2)\z(3)-W(5,2)=\frac{21}{16}\z(2)\z(3)-\frac{31}{16}\z(5),
\end{align*}
where we used the formula (\ref{b24}) and the identity
\[W(5,2)=\frac{7}{16}\z(2)\z(3)+\frac{31}{16}\z(5).\]

\begin{thm} For positive integers $m$ and $p$,
\begin{align}\label{c4}
2\su \frac{T_n(\{1\}_{2m-2})T_n(\{1\}_{2p-2})}{2n-1}(-1)^n=(-1)^{m-1} \sum_{j=1}^p W(2j+2m-3,2m-1)Z(j,p).
\end{align}
\end{thm}
\pf  Setting
\begin{align*}
&A_{j,p}:=2{\bar \z}(2p-2j),\quad B_j:=2\su \frac{T_n(\{1\}_{2m-2})T_n(\{1\}_{2j-2})}{2n-1}(-1)^n
\end{align*}
and
\begin{align*}
C_p:=(-1)^{m-1}W(2p+2m-3,2m-1)
\end{align*}
in Lemma \ref{lemc1} and using (\ref{b23}) yields the desired formula.\hfill$\square$

Letting $m=p=2$ in (\ref{c4}), we have
\begin{align*}
2\su \frac{T^2_n(1,1)}{2n-1}(-1)^n=-W(5,3)+2{\bar \z}(2)\T(1,1,1).
\end{align*}

Applying (\ref{b21}) and (\ref{b22}), we know that for positive integers $m$ and $p$, the sums
\[\su \frac{T_n(\{1\}_{2m-1})S_n(\{1\}_{2p-2})}{n}(-1)^n\quad \text{and}\quad \su \frac{T_n(\{1\}_{2m-2})S_n(\{1\}_{2p-1})}{2n-1}(-1)^n\]
can also be evaluated by weighted sums $W(k,r)$ and $Z(j,p)$.
\\
{\bf Proofs of Theorems \ref{thm1} and \ref{thm2}}. Changing $(m,p)$ to $(p,m)$ in (\ref{c3}) and (\ref{c4}), and using the duality of series on the left hand sides, then
\begin{align*}
\su \frac{T_n(\{1\}_{2m-1})T_n(\{1\}_{2p-1})}{n}(-1)^n&=(-1)^m \sum_{j=1}^p W(2j+2m-1,2m)Z(j,p)\\
&=(-1)^p \sum_{j=1}^m W(2j+2p-1,2p)Z(j,m)
\end{align*}
and
\begin{align*}
2\su \frac{T_n(\{1\}_{2m-2})T_n(\{1\}_{2p-2})}{2n-1}(-1)^n&=(-1)^{m-1} \sum_{j=1}^p W(2j+2m-3,2m-1)Z(j,p)\\
&=(-1)^{p-1} \sum_{j=1}^m W(2j+2p-3,2p-1)Z(j,m).
\end{align*}
Thus, we prove Theorems \ref{thm1} and \ref{thm2}. \hfill$\square$

We get some cases. Setting $m=1$ in (\ref{a9}) and (\ref{a10}) yield
\begin{align}\label{c5}
W(2p+1,2p)=(-1)^{p-1} \sum_{j=1}^p W(2j+1,2)Z(j,p)
\end{align}
and
\begin{align}
\T(\{1\}_{2p-1})=(-1)^{p-1}\sum_{j=1}^p \T(2j-1)Z(j,p).
\end{align}

Now, we give the explicit expression of weighted sums $W(2k+1,2)$ via $t$-values and alternating zeta values.
\begin{thm}\label{thmc4} For positive integer $k$,
\begin{align}\label{c7}
W(2k+1,2)=\frac{1}{2^{2k-1}} \sum_{j=1}^k \t(2j+1){\bar \z}(2k-2j),
\end{align}
where $\t(k):=2^k t(k)=(2^k-1)\z(k)$ if $k>1$ and $\t(1):=0$ if $k=1$.
\end{thm}
\pf From \cite[Eq.(3.10)]{WX2020}, Wang and the author gave the result
\begin{align}\label{c8}
&(1-(-1)^{p+q})\su \frac{h^{(p)}_n}{n^q}(-1)^{n-1}\nonumber\\
=&-(-1)^p(1+(-1)^q)\t(p){\bar \z}(q)+(-1)^p\binom{p+q-1}{p-1}\t(p+q)\nonumber\\
&-(-1)^p\sum_{k=0}^{p-1} ((-1)^k+1)\binom{p+q-k-2}{q-1}{\bar t}(k+1){\bar t}(p+q-k-1)\nonumber\\
&+2(-1)^p \sum_{j=1}^{[q/2]} \binom{p+q-2j-1}{p-1}{\bar \z}(2j)\t(p+q-2j).
\end{align}
According to definition, we have
\begin{align*}
W(2k+1,2)=\sum_{k_1+k_2=2k+1,\atop k_1,k_2\geq 1} \frac{1}{2^{k_1+k_2-2}} \su \frac{h^{(k_1)}_n}{n^{k_2}}(-1)^{n-1}.
\end{align*}
Thus, applying (\ref{c8}) with $p=k_1$ and $q=k_2$, by an elementary calculation, we complete this proof. \hfill$\square$

Hence, using (\ref{c5}) and (\ref{c7}), we obtain the following description.
\begin{cor} For positive integer $p$, the weighted sums $W(2p+1,2p)$ can be expressed in terms of (alternating) Riemann zeta values.
\end{cor}

As two examples, for $p=1$ and $2$, we have
\begin{align*}
W(3,2)=\frac{7}{4}\z(3)\quad \text{and} \quad W(5,4)=\frac{21}{16}\z(2)\z(3)-\frac{31}{16}\z(5).
\end{align*}

\section{Proof of Theorem \ref{thm3}}
In this section we prove Theorem \ref{thm3}. First, we need to give a lemma.
\begin{lem}\label{lem4.1}(\cite[Lemma 2.5]{X2020}) If $f_i\ (i=1,\ldots,m)$ are integrable real functions, the following identity holds:
\begin{align}
 &g\left( {{f_1},{f_2}, \cdots ,{f_m}} \right) + {\left( { - 1} \right)^m}g\left( {{f_m},{f_{m - 1}}, \cdots ,{f_1}} \right) \nonumber\\
 & = \sum\limits_{i = 1}^{m - 1} {{{\left( { - 1} \right)}^{i - 1}}g\left( {{f_i},{f_{i - 1}}, \cdots ,{f_1}} \right)} g\left( {{f_{i + 1}},{f_{i + 2}} \cdots ,{f_m}} \right),
\end{align}
where $g\left( {{f_1},{f_2}, \cdots ,{f_m}} \right)$ is defined by
\[g\left( {{f_1},{f_2}, \cdots ,{f_m}} \right): = \int\limits_{0 < {t_m} <  \cdots <t_2 < {t_1} < 1} {{f_1}\left( {{t_1}} \right){f_2}\left( {{t_2}} \right) \cdots {f_m}\left( {{t_m}} \right)d{t_1}d{t_2} \cdots d{t_m}} .\]
\end{lem}

On the one hand, in (\ref{b5}), replacing $k_r$ by $k_r+m-1$ and letting $z=i$, then summing, we find that
\begin{align}
&W_m(k+r-1,r)=\sum_{{\bf k}\in I(k+r-1,r)} \binom{k_r+m-2}{m-1} \T(k_1,\ldots,k_{r-1},k_r+m-1)\nonumber\\
&=\frac{(-1)^{k+r+m}}{(k-1)!(m-1)!}2^r \int_{D_r(1)} \frac{\log^{m-1}\left(\frac{1-t_1}{1+t_1} \right)\log^{k-1}\left(\frac{1-t_r}{1+t_r} \right)}{(1+t_1^2)\cdots (1+t_r^2)}dt_1\cdots dt_r.
\end{align}
On the other hand, we deduce the identity
\begin{align}\label{d3}
\int_{D_r(1)} \frac{\log^{p}\left(\frac{1-t_1}{1+t_1} \right)}{(1+t_1^2)\cdots (1+t_r^2)}dt_1\cdots dt_r
=\frac{(-1)^{p+r}p!}{2^r} T(\{1\}_{p-1},{\bar 1},\{1\}_{r-1},{\bar 1})
\end{align}
by using (\ref{b1}), (\ref{b14}) and expanding the $(1+t_j^2)^{-1}$ in geometric series with a direct calculation.

Then using Lemma \ref{lem4.1}, by an elementary calculation, we prove the Theorem \ref{thm3}.\hfill$\square$

Letting $m=1$ in Theorem \ref{thm3}, we can get the following a corollary.
\begin{cor} For positive integers $k$ and $r$,
\begin{align}\label{d4}
W(k+r-1,r)+(-1)^r\T(\{1\}_{r-1},k)=\sum_{j=1}^{r-1}(-1)^{j-1}T(\{1\}_{k-2},{\bar 1},\{1\}_{j-1},{\bar 1})\T(\{1\}_{r-j}).
\end{align}
\end{cor}

Setting $r=2$ in (\ref{d4}) yields
\begin{align}\label{d5}
W(k+1,2)+\T(1,k)=T(\{1\}_{k-2},{\bar 1},{\bar 1})\T(1).
\end{align}
where $\T(1)=-{\pi}/{2}$ from (\ref{b15}).

Next, we give a duality theorem of alternating MTVs and give some weighted sum formulas involving alternating MTVs.
\begin{thm} For integers $p\geq 0$ and $r\geq 1$,
\begin{align}\label{d6}
T(\{1\}_{p-1},{\bar 1},\{1\}_{r-1},{\bar 1})=\T(\{1\}_{r-1},p+1),
\end{align}
where $(\{1\}_{-1},{\bar 1}):=\emptyset$.
\end{thm}
\pf In (\ref{d3}), applying the changes of variables $t_j\rightarrow \frac{1-t_{r+1-j}}{1+t_{r+1-j}}$ to (\ref{b4}) gives
\begin{align*}
\frac{(-1)^{p+r}p!}{2^r} T(\{1\}_{p-1},{\bar 1},\{1\}_{r-1},{\bar 1})&=\int_{D_r(1)} \frac{\log^p(t_r)}{(1+t_1^2)\cdots (1+t_r^2)}dt_1\cdots dt_r\\
&=\frac{(-1)^{p+r}p!}{2^r} \T(\{1\}_{r-1},p+1),
\end{align*}
where in the last step, we expanded the $(1+t_j^2)^{-1}$ in geometric series and used the following well-known identity
$$\int_{0}^1 t^{n-1} \log^p(t)dt=p!\frac{(-1)^p}{n^{p+1}}.$$
Thus, we complete this proof. \hfill$\square$

Letting $r=1$ in (\ref{d6}) yields
\begin{align}\label{d7}
T(\{1\}_{p-1},{\bar 1},{\bar 1})=-\frac{1}{2^p}{\bar t}(p+1).
\end{align}
In particular, if $k=2$, then we have
\[T({\bar 1},{\bar 1})=-\frac{1}{2}{\bar t}(2)=-2G,\]
where $G:=\sum_{n=1}^\infty (-1)^{n-1}/(2n-1)^2$ is Catalan's constant (see \cite[Eq. (7.4.4.183a)]{K1997}).
Hence, substituting (\ref{d7}) into (\ref{d5}) with $p=k-1$, we get the following corollary.
\begin{cor} For positive integer $k$,
\begin{align}\label{d8}
W(k+1,2)=\frac{\pi}{2^k}{\bar t}(k)-\T(1,k).
\end{align}
\end{cor}

\begin{re} The formula (\ref{d8}) can be obtained by Theorem \ref{thmb3} with $r=2$.
\end{re}

Further, applying (\ref{d6}) with $p=k-1$, the formula (\ref{a11}) can be rewritten in the form
\begin{align}\label{d9}
W_m(k+r-1,r)+(-1)^rW_k(m+r-1,r)=\sum_{j=1}^{r-1}(-1)^{j-1}\T(\{1\}_{j-1},k)\T(\{1\}_{r-j-1},m)
\end{align}
and
the formula (\ref{d4}) can be rewritten in the form
\begin{align}\label{d10}
W(k+r-1,r)+(-1)^r\T(\{1\}_{r-1},k)=\sum_{j=1}^{r-1}(-1)^{j-1}\T(\{1\}_{j-1},k)\T(\{1\}_{r-j}).
\end{align}
Clearly, (\ref{d10}) can also be immediately obtained from (\ref{d9}) with $m=1$. Moreover, setting $r=2$ in (\ref{d10}) yields (\ref{d8}).

Now, we end this section by a theorem.
\begin{thm}\label{thmd6} For positive integers $m,k$ and $r$,
\begin{align}\label{d11}
&\sum_{k_1+\cdots+k_r=k+r-1,\atop k_1,\ldots,k_r\geq 1} \binom{k_r+m-2}{m-1}{\rm A}(k_1,\ldots,k_{r-1},k_r+m-1;z)\nonumber\\
&+(-1)^r \sum_{k_1+\cdots+k_r=m+r-1,\atop k_1,\ldots,k_r\geq 1} \binom{k_r+k-2}{k-1}{\rm A}(k_1,\ldots,k_{r-1},k_r+k-1;z)\nonumber\\
&=\sum_{j=1}^{r-1} (-1)^{j-1} {\rm A}(\{1\}_{r-1-j},m;z){\rm A}(\{1\}_{j-1},k;z).
\end{align}
\end{thm}
\pf Replacing $k_r$ by $k_r+m-1$ in (\ref{b5}) and summing, by an elementary calculation, we have
\begin{align}\label{d12}
&\sum_{k_1+\cdots+k_r=k+r-1,\atop k_1,\ldots,k_r\geq 1} \binom{k_r+m-2}{m-1} {\rm A}(k_1,\ldots,k_{r-1},k_r+m-1;z)\nonumber\\
&=\frac{(-1)^{k+m}}{(m-1)!(k-1)!} 4^r z^r \int_{D_r(1)} \frac{\log^{m-1}\left(\frac{1-t_1}{1+t_1} \right)\log^{k-1}\left(\frac{1-t_r}{1+t_r} \right)dt_1\cdots dt_r}{[(1+t_1)^2-z^2(1-t_1)^2]\cdots [(1+t_r)^2-z^2(1-t_r)^2] }.
\end{align}
Letting $k=1$ yields
\begin{align*}
{\rm A}(\{1\}_{r-1},m;z)=\frac{(-1)^{m-1}}{(m-1)!} 4^r z^r \int_{D_r(1)} \frac{\log^{m-1}\left(\frac{1-t_1}{1+t_1} \right)dt_1\cdots dt_r}{[(1+t_1)^2-z^2(1-t_1)^2]\cdots [(1+t_r)^2-z^2(1-t_r)^2] }.
\end{align*}
Finally, using Lemma \ref{lem4.1}, by a direct calculation, we may easily deduce the desired result. \hfill$\square$

It is clear that setting $z=-1$ in (\ref{d11}) gives the identity (\ref{d9}).
\\[5mm]
{\bf Acknowledgments.}  The author expresses his deep gratitude to Professors Masanobu Kaneko and Jianqiang Zhao for valuable discussions and comments.

 \end{document}